\documentclass{article}
\usepackage[utf8]{inputenc}
\usepackage{graphicx}
\graphicspath{ {images/} }
\usepackage{caption}
\usepackage{subcaption}
\usepackage{float}
\usepackage[width=150mm,top=35mm,bottom=25mm,bindingoffset=6mm]{geometry}
\usepackage{fancyhdr}
\usepackage{setspace}
\usepackage{amssymb, amsmath}

\usepackage{amssymb,amsmath,amsthm}
\usepackage{esint}
\usepackage{mathrsfs}
\usepackage{mathtools}
\usepackage{comment}
\usepackage{bbm,dsfont}

\usepackage{graphicx}
\usepackage[utf8]{inputenc}
\usepackage[T1]{fontenc}
\usepackage[hidelinks]{hyperref}

\usepackage{comment}

\usepackage[numbers]{natbib}
\newtheorem{definition}{Definition}

\newtheorem{theorem}{Theorem}
\newtheorem{corollary}{Corollary}
\newtheorem{proposition}{Proposition}
\usepackage{graphicx} 
\title{Generalized Variance Inequalities for Barycenters in CAT(0) and CAT(1) Spaces} 
\author{Sebastian Gietl}
\date{}

\begin{document}
\maketitle
\begin{abstract}
We prove generalized versions of the Variance Inequality known for barycenters in CAT(0) spaces, inspired by an analogous result for $p$-uniformly convex Banach spaces. Our generalizations apply to balls of sufficiently small radius in complete CAT(1) spaces and to exponents $p \geq 2$ in the $\operatorname{CAT}(0)$ setting. Building on a result of Eskenazis, Mendel, and Naor, we establish sharp metric cotype for all $p \geq 2$ in $\mathrm{CAT}(0)$ spaces, extending the previously known case $p=2$. In addition, based on their work, we derive martingale inequalities for nonlinear martingales taking values in complete $\mathrm{CAT}(0)$ space and balls of sufficiently small radius in complete CAT(1) spaces.
\end{abstract}

\section{Introduction}\label{intro}
Let \((X, d_X)\) be a complete \(\operatorname{CAT}(0)\) space. The well-known Variance Inequality characterizes such metric spaces, for any square-integrable Borel probability measure \(\mu\) on \(X\), there exists a point \(z_\mu \in X\) such that
\[
\int d(z, x)^2 \, d\mu(x) \geq d\left(z, z_{\mu}\right)^2 + \int d\left(z_{\mu}, x\right)^2 \, d\mu(x) \quad \forall z \in X.
\]
See Theorem 4.9 in \cite{PMonNPC}. In fact, the point \( z_\mu \) can be taken to be the barycenter \( b(\mu) \) of \( \mu \) (see below for the definition). In this note, we aim to extend this inequality to the \(\operatorname{CAT}(1)\) setting and, in the case of nonpositive curvature, to more general exponents \( p \geq 2 \). These generalized inequalities mirror a result established in \cite{MCandLip}, known to hold in \( p \)-uniformly convex Banach spaces (see, for instance, Chapter 3 of \cite{MCandLip} for a definition).
\begin{theorem}[Generalized Variance Inequality for \(\operatorname{CAT}(1)\) spaces]\label{VarInCAT1}
Let \((X, d_X)\) be a complete \(\operatorname{CAT}(1)\) space, and let \(\mu\) be a square-integrable Borel probability measure supported on a ball \(B_r(o) \subset X\) of radius \(0 < r < \pi/2\), with midpoint \(o \in X\). Then,
\[
\int d_X(o, x)^2 \, d\mu(x) \geq d_X(o, b(\mu))^2 + \frac{k_r}{2} \int d_X(b(\mu), x)^2 \, d\mu(x),
\]
where \(k_r = 2r \tan\left(\frac{\pi}{2} - r\right)\).
\end{theorem}
Note that \( \lim_{r \to 0} k_r = 2 \), so in the limit we recover the \(\operatorname{CAT}(0)\) constant. Additionally, compared to the \(\operatorname{CAT}(0)\) version above, the arbitrary point \( z \) is replaced by the midpoint of the ball. It is possible to prove a version of Theorem \ref{VarInCAT1} with an arbitrary point \( z \), however, this requires assuming that the diameter, rather than just the radius, of the space is strictly less than \( \pi/2 \).
\begin{theorem}[Generalized Variance Inequality for \(\operatorname{CAT}(0)\) spaces]\label{VarInCAT(0)}
Let \(p \in [2, \infty)\), and let \((X, d_X)\) be a complete \(\operatorname{CAT}(0)\) space. Then, for any \(p\)-integrable Borel probability measure \(\mu\) on \(X\), we have
\[
\int d_X(z, x)^p \, d\mu(x) \geq d_X(z, b(\mu))^p + \frac{2^{p-2}}{2^{p-1} - 1} \cdot \frac{k_p}{2} \int d_X(b(\mu), x)^p \, d\mu(x) \quad \forall z \in X,
\]
where $k_p=2$ if $p=2$ and $k_p=\left(\frac{(p-1)^2}{10 p^2}\right)^{p-1}$ if $p>2$.
\end{theorem}

We will establish a more general version of Theorems~\ref{VarInCAT1} and~\ref{VarInCAT(0)} for metric spaces in which the \(p\)-th power of the distance is uniformly convex and which admit a suitable barycenter map. A different generalization of the Variance Inequality was given in \cite{JenOnConSp}, but it does not follow the Banach space analogy and therefore does not lead to the same applications as the one developed here.

As an application of Theorems~\ref{VarInCAT1} and \ref{VarInCAT(0)}, and using a result from \cite{LipForBary}, we derive nonlinear versions of a classical martingale inequality due to Pisier, originally established in \cite{MartInUniB} for uniformly convex Banach spaces. Furthermore, combining our results with a theorem of Eskenazis, Mendel, and Naor from \cite{CATCoarU}, we obtain that \(\operatorname{CAT}(0)\) spaces possess sharp metric cotype \( p \) for all \( p \geq 2 \), a nonlinear analogue of the Banach space notion of Rademacher cotype introduced in \cite{MC}. In \cite{CATCoarU}, this was originally shown for the case \( p = 2 \).

\section{Preliminaries}\label{AlexSp}
\subsection{CAT(0) and CAT(1) Spaces}
 We briefly recall some facts about $\operatorname{CAT}(1)$ spaces. Such spaces were treated systematically in \cite{MSofNPC}. Let \((X, d_X)\) be a metric space and \(I \subset \mathbb{R}\) be an interval. A curve \(\gamma: I \rightarrow X\) is called a geodesic if  \(d_X(\gamma(s), \gamma(t)) =  |s - t|\) for all \(s, t \in I\). For \(x, y \in X\), a curve \(\gamma: [0, 1] \rightarrow X\) is said to be a geodesic connecting \(x\) and \(y\) when it is a geodesic and \(\gamma(0) = x\) and \(\gamma(1) = y\). The space \((X, d_X)\) is said to be geodesic if, for every pair of points \(x, y \in X\), there exists a geodesic connecting \(x\) and \(y\). If $(X,d_X)$ is geodesic, a function \(f:X\rightarrow\mathbb{R}\) is said to be convex if for each geodesic \(\gamma:[0,1]\rightarrow X \) the composition \(f \circ \gamma:[0,1] \rightarrow \mathbb{R}\) is convex. 
 
Consider \(\kappa \in \mathbb{R}\), we define \(\mathbb{M}^2(\kappa)\) as the 2-dimensional space form of constant curvature \(\kappa\). Then the diameter of \(\mathbb{M}^2(\kappa)\), denoted by \(D_\kappa\), is given by:
\begin{equation*}
    D_\kappa =
\begin{cases}
\infty & \text{if } \kappa \leq 0, \\
\frac{\pi}{\sqrt{\kappa}} & \text{if } \kappa > 0.
\end{cases}
\end{equation*}
We say that \((X, d_X)\) is $D_{\kappa}$-geodesic if, for every pair of points \(x, y \in X\) with $d(x,y)<D_{\kappa}$, there exists a geodesic connecting \(x\) and \(y\).  The image of a geodesic connecting \(x\) and \(y\) is called a geodesic segment and denoted by \([x,y]\). A geodesic triangle \(\Delta\) in \(X\) consists of three points \(x, y, z \in X\), its vertices, and a choice of three geodesic segments \([x, y],[y, z],[z, x]\) joining them, its sides. A triangle \(\bar{\Delta}\) in \(\mathbb{M}^2(\kappa)\) with vertices \(\bar{x}, \bar{y}, \bar{z} \in \mathbb{M}^2(\kappa)\) is called a comparison triangle for \(\Delta\) if \(d_{\mathbb{M}^2(\kappa)}(\bar{x}, \bar{y})=d_X(x, y), d_{\mathbb{M}^2(\kappa)}(\bar{y}, \bar{z})=d_X(y, z)\) and \(d_{\mathbb{M}^2(\kappa)}(\bar{z}, \bar{x})=d_X(z, x)\). Such a triangle \(\bar{\Delta} \subseteq \mathbb{M}^2(\kappa)\) always exists if the perimeter \(d_X(x, y) + d_X(y, z) + d_X(z, x)\) of \(\Delta\) is less than \(2 D_\kappa\), and it is unique up to isometry. The space \((X, d_X)\) is said to be a $\operatorname{CAT}(\kappa)$ space if it is $D_{\kappa}$-geodesic and for all geodesic triangles in \(X\) with perimeter strictly less than \(2D_\kappa\), it holds that for all \(x, y \in \Delta\) and their corresponding points \(\bar{x}, \bar{y} \in \bar{\Delta}\), the distance satisfies
\begin{equation}\label{distComp}
    d_X\left(x, y\right) \leq d_{\mathbb{M}^2(\kappa)}\left(\bar{x}, \bar{y}\right).
\end{equation}
Note that for \(\kappa_1, \kappa_2 \in \mathbb{R}\) with \(\kappa_2 \geq \kappa_1\), any \(\operatorname{CAT}(\kappa_1)\) space is also a \(\operatorname{CAT}(\kappa_2)\) space. Moreover, by rescaling the metric, it suffices to consider the cases \(\kappa = -1\), \(\kappa = 0\), and \(\kappa = 1\). 

If \((X, d_X)\) is a \(\operatorname{CAT}(0)\) space, then equivalently, for all points \(x, y, z \in X\) and any geodesic \(\gamma: [0,1] \rightarrow X\) connecting \(x\) and \(y\), we have
\begin{equation}\label{con0}
d_X\left(z, \gamma\left(\tfrac{1}{2}\right)\right)^2 \leq \tfrac{1}{2}\, d_X\left(z, x\right)^2 + \tfrac{1}{2}\, d_X\left(z, y\right)^2 - \tfrac{1}{4}\, d_X\left(x, y\right)^2
\end{equation}
holds. A generalization for exponents \(p \geq 2\) takes the form
\begin{equation}\label{conp}
d_X\left(z, \gamma\left(\tfrac{1}{2}\right)\right)^p \leq \tfrac{1}{2}\, d_X\left(z, x\right)^p + \tfrac{1}{2}\, d_X\left(z, y\right)^p - \tfrac{k_p}{2}  \tfrac{1}{4}\, d_X\left(x, y\right)^p,
\end{equation}
where \( k_p = \left( \frac{(p - 1)^2}{10p^2} \right)^{p - 1} \). This can be deduced from \cite[Corollary 6.4]{SpeCalandSup}, together with the fact that \( \mathbb{R}^2 \) embeds isometrically into any \( L_p \)-space. If \(X = B_r(o)\) is a ball of radius \(0 < r < \pi/2\) centered at a midpoint \(o\) in a \(\operatorname{CAT}(1)\) space, then a version of a result originally shown by Ohta in \cite[Proposition 3.1]{ConOfms}, and made explicit in \cite[Lemma 5]{BaryOnCat(1)}, states the following, for all \(x, y \in X\) and all geodesics \(\gamma : [0,1] \to X\) connecting \(x\) and \(y\), we have
\begin{equation}\label{conk}
d_X\left(o, \gamma\left(\tfrac{1}{2}\right)\right)^2 \leq \tfrac{1}{2} d_X\left(o, x\right)^2 + \tfrac{1}{2} d_X\left(o, y\right)^2 - \tfrac{k_r}{2} \tfrac{1}{4} d_X\left(x, y\right)^2.
\end{equation}
where \(k_{r} = 2r \tan \left(\frac{\pi}{2} - r\right)\). Ohta's original result establishes Inequality~\ref{conk} with \( o \) replaced by an arbitrary point in \( X \), but assumes that \( \operatorname{diam}(X) < \pi/2 \).

\subsection{Barycenters}
Barycenters have been studied, for example, in \cite{PMonNPC} for \(\operatorname{CAT}(0)\) spaces and in \cite{BaryOnCat(1)} for \(\operatorname{CAT}(1)\) spaces. Let \((X, d_X)\) be a metric space and \(p \in [1, \infty)\). We denote by \(\mathcal{P}_p(X)\) the set of all Borel probability measures on \(X\) with finite \(p\)-th moment; that is, for all \(\mu \in \mathcal{P}_p(X)\), there exists some \(z \in X\) such that
\[
\int d_X(z, x)^p \, d\mu(x) < \infty,
\]
which then holds for all \(z \in X\). We denote by \(\mathcal{P}_\infty(X)\) the set of all Borel probability measures on \(X\) with finite support. 

 Let \( (X, d_X) \) be a complete \( \operatorname{CAT}(0) \) space or a ball $B_r $ of radius $0<r<\pi/ 2$ in a complete $\operatorname{CAT}(1)$ space and $\mu \in \mathcal{P}_2(X)$. Then the functional
\[
z \mapsto \int d_X(z, x)^2 \, d\mu(x)
\]
admits a unique minimizer \(b(\mu)\), called the barycenter of \(\mu\). The barycenter satisfies a Jensen inequality. That is, for any lower semicontinuous convex function \(\varphi: X \to \mathbb{R}\) that is \(\mu\)-integrable, we have
\[
\varphi(b(\mu)) \leq \int \varphi(x) \, d\mu(x).
\]
The \(\operatorname{CAT}(0)\) case is proved, for example, in \cite[Theorem 6.2]{PMonNPC}, while the \(\operatorname{CAT}(1)\) case is addressed in \cite[Theorem 25]{BaryOnCat(1)}.

\section{The Generalized Variance Inequality}\label{bary}
The concept of $p$-uniformly convex metric spaces was studied by Naor and Silberman in \cite{PIandWG}, and by Ohta in \cite{ConOfms}.
\begin{definition}[$p$-uniform convexity]
Let \((X, d_X)\) be a geodesic space and let \(p \in [2, \infty)\). The space \((X, d_X)\) is said to be \(p\)-uniformly convex with respect to $z\in X$ with constant \(k > 0\) if, for every \(x, y\in X\) and every geodesic \(\gamma: [0,1] \rightarrow X\) connecting \(x\) and \(y\), we have
\[
d_X\left(z, \gamma\left(\tfrac{1}{2}\right)\right)^p \leq \tfrac{1}{2} d_X(z, x)^p + \tfrac{1}{2} d_X(z, y)^p - \tfrac{k}{2}  \tfrac{1}{4} d_X(x, y)^p.
\]
We say that \((X, d_X)\) is \(p\)-uniformly convex with constant \(k > 0\) if, for all \(z \in X\), the space \((X, d_X)\) is \(p\)-uniformly convex with respect to \(z\) with constant \(k\).
\end{definition}
We will make use of the fact that geodesics are unique in \(p\)-uniformly convex metric spaces; see, for example, \cite[Lemma 2.3]{ConOfms}. Our goal is to prove a generalized variance inequality for \(p\)-uniformly convex metric spaces that admit a suitable barycenter map. To this end, we introduce the following definition.  
\begin{definition}[Barycenter Map]
Let \((X, d_X)\) be a metric space and \(p \in [1, \infty]\). A map \(\mathfrak{B}: \mathcal{P}_p(X) \rightarrow X\) is called a \(p\)-barycenter map if for all \(x \in X\), we have \(\mathfrak{B}(\delta_x) = x\).
\begin{enumerate}
\item A \(p\)-barycenter map \(\mathfrak{B}\) is called a \(p\)-convex mean map if, for every \(\mu \in \mathcal{P}_p(X)\) and every lower semicontinuous, \(\mu\)-integrable, and convex function \(\varphi: X \rightarrow \mathbb{R}\), the following Jensen-type inequality holds:
\[
    \varphi(\mathfrak{B}(\mu)) \leq \int \varphi \, d\mu.
\]
\item Assume that \((X, d_X)\) is uniquely geodesic. For \(\mu \in \mathcal{P}_p(X)\), we denote by \(m_{\mu} : X \to X\) the map that sends a point \(x \in X\) to the unique geodesic midpoint between \(x\) and \(\mathfrak{B}(\mu)\). The $p$-barycenter map \(\mathfrak{B}\) is said to be invariant under midpoint contractions if, for every \(\mu \in \mathcal{P}_p(X)\), we have
\[
\mathfrak{B}(\mu) = \mathfrak{B}((m_{\mu})_*\mu).
\]
\end{enumerate}
\end{definition}
We now turn our attention to the generalized variance inequality. Our proof is a metric version of Lemma 3.1 in \cite{MCandLip}, where the result was demonstrated for $2$-uniformly convex Banach spaces, in which the barycenter map is defined as $\mu \mapsto \int x \, d\mu(x)$. A generalized version for $p$-uniformly convex Banach spaces with $p \in [2, \infty)$ appears in Lemma 6.5 of \cite{SpeCalandSup}.
\begin{theorem}[Generalized Variance Inequality for $p$-uniformly convex spaces]\label{VarIn}
Let \(p \in [2, \infty)\), let \((X, d_X)\) be a uniquely geodesic space, and let \(z \in X\). Assume that $(X, d_X)$ is $p$-uniformly convex with respect to $z$ with constant $k > 0$, which admits a $p$-convex mean map $\mathfrak{B} \colon \mathcal{P}_p(X) \to X$ that is invariant under midpoint contractions. Then, for any $\mu \in \mathcal{P}_p(X)$, we have
\begin{equation}\label{VarInEq}
\int d_X(z, x)^p \, d\mu(x) \geq d_X(z, \mathfrak{B}(\mu))^p + \frac{2^{p-2}}{2^{p-1} - 1} \cdot \frac{k}{2} \int d_X(\mathfrak{B}(\mu), x)^p \, d\mu(x).
\end{equation}
In particular, if \((X, d_X)\) is a \(p\)-uniformly convex metric space with constant \(k > 0\) that admits a barycenter map as above, then Equation~\eqref{VarInEq} holds for every \(z \in X\).
\end{theorem}

\begin{proof}
 Because $\mathfrak{B}$ is a $p$-convex mean map, the following infimum exists,
\begin{equation*}
    \theta := \inf \left\{\frac{\int d(z, x)^p d \mu(x) - d(z, \mathfrak{B}(\mu))^p}{\int d(x, \mathfrak{B}(\mu))^p d \mu(x)} : \mu \in \mathcal{P}_p(X) \text{ s.t. } \mu \neq \delta_x \quad \forall x \in X \right\}.
\end{equation*}
So, for any $\phi > \theta$, there exists $\mu_0 \in \mathcal{P}_p(X)$ such that
\begin{equation}\label{1.1}
    \phi \int d(x, \mathfrak{B}(\mu_0))^p d \mu_0(x) > \int d\left(z, x\right)^p d \mu_0(x) - d(z, \mathfrak{B}(\mu_0))^p.
\end{equation}
Also, for every \(x \in X\), we have
\begin{equation}\label{1.2}
    d(z, x)^p \geq 2d\left(z, m_{\mu_0}(x)\right)^p + \frac{k}{4} d\left(x, \mathfrak{B}(\mu_0)\right)^p - d\left(z, \mathfrak{B}(\mu_0)\right)^p.
\end{equation}
 Combining (\ref{1.1}) and (\ref{1.2}) together with the fact that $\mathfrak{B}$ is invariant under midpoint contractions yields
\begin{equation*}
\begin{aligned}
    \phi \int &d\left(x, \mathfrak{B}(\mu_0)\right)^p d \mu_0(x) \\
    &> 2\left(\int d\left(z, m_{\mu_0}(x)\right)^p d \mu_0(x) - d(z, \mathfrak{B}(\mu_0))^p\right) + \frac{k}{4} \int d\left(x, \mathfrak{B}(\mu_0)\right)^p d\mu_0(x) \\
    &= 2\left(\int d\left(z, m_{\mu_0}(x)\right)^p d \mu_0(x) - d\left(z, \mathfrak{B}\left(\left(m_{\mu_0}\right)_* \mu_0\right)\right)^p\right) + \frac{k}{4} \int d\left(x, \mathfrak{B}(\mu_0)\right)^p d\mu_0(x) \\
    &\geq 2 \theta \int d\left(m_{\mu_0}(x), \mathfrak{B}(\mu_0)\right)^p d \mu_0(x) + \frac{k}{4} \int d\left(x, \mathfrak{B}(\mu_0)\right)^p d \mu_0(x) \\
    &= \frac{\theta}{2^{p-1}} \int d\left(x, \mathfrak{B}(\mu_0)\right)^p d \mu_0(x) + \frac{k}{4} \int d\left(x, \mathfrak{B}(\mu_0)\right)^p d \mu_0(x).
\end{aligned}
\end{equation*}
Therefore, for any $\phi > \theta$, it holds that $\phi > \left(\frac{\theta}{2^{p-1}} + \frac{k}{4}\right)$, which implies $\theta \geq \frac{2^{p-2}}{2^{p-1}-1}\frac{k}{2}$.
\end{proof}
Now we show that the barycenter map defined on CAT(1) spaces is invariant under midpoint contractions
\begin{proposition}\label{PushLem}
Let \( (X, d_X) \) be a complete \( \operatorname{CAT}(0) \) space or a ball $B_r $ of radius $0<r<\pi/ 2$ in a complete $\operatorname{CAT}(1)$ space. Then the barycenter map $b:\mathcal{P}_2(X)\rightarrow X$ is invariant under midpoint contractions.
\end{proposition}
\begin{proof}
Fix $z \in X$ and $\mu \in \mathcal{P}_2(X)$. For every \(x \in X\), we have
\begin{equation*}
    \begin{aligned}
        & \int d_X(z, x)^2 \, d\mu(x) \geq \int d_X(b(\mu), x)^2 \, d\mu(x) = \int 4\, d_X\left(b(\mu), m_{\mu}(x)\right)^2 \, d\mu(x) \\
        & = 2 \int d_X\left(b(\mu), m_{\mu}(x)\right)^2 \, d\mu(x) + 2 \int d_X\left(m_{\mu}(x), x\right)^2 \, d\mu(x) \\
        & \geq 2 \int d_X\left(b\left(\left(m_\mu\right)_* \mu\right), m_{\mu}(x)\right)^2 \, d\mu(x) + 2 \int d_X\left(m_{\mu}(x), x\right)^2 \, d\mu(x) \\
        & = 4 \int \left[ \tfrac{1}{2} d_X\left(b\left(\left(m_\mu\right)_* \mu\right), m_{\mu}(x)\right)^2 + \tfrac{1}{2} d_X\left(m_{\mu}(x), x\right)^2 \right] d\mu(x) \\
        & \geq 4 \int \tfrac{1}{4} d_X\left(b\left(\left(m_\mu\right)_* \mu\right), x\right)^2 \, d\mu(x) = \int d_X\left(b\left(\left(m_\mu\right)_* \mu\right), x\right)^2 \, d\mu(x),
    \end{aligned}
\end{equation*}
where the last inequality is a consequence of the triangle inequality. By the uniqueness of barycenters, the claim follows.
\end{proof}
Theorems~\ref{VarInCAT1} and~\ref{VarInCAT(0)} follow from Proposition~\ref{PushLem} and Theorem~\ref{VarIn}, combined with the appropriate geometric properties in each curvature setting. For the \(\operatorname{CAT}(1)\) case, we use property~\eqref{conk} and the fact that geodesics in balls of radius \(r < \pi/2\) are unique (see \cite[II.1.4]{MSofNPC}). For the \(\operatorname{CAT}(0)\) case, we rely on property~\eqref{conp}.

\section{Applications}\label{NC}
In this section, we formulate applications of the generalized variance inequalities. In particular, we examine their use in establishing metric cotype, a nonlinear analogue of Rademacher cotype introduced in \cite{MC}. To define metric cotype, we begin with some notation. Let \(\mathbb{Z}_{2m}\) denote the set of integers modulo \(2m\), and note that additions appearing in the definition below are performed modulo \(2m\). Furthermore, let \(e_1 = (1, 0, \ldots, 0), \ldots, e_n = (0, \ldots, 0, 1)\) denote the standard basis vectors of \(\mathbb{Z}_{2m}^n\).

\begin{definition}[Metric Cotype]\label{MC}
A metric space \((X, d_X)\) is said to have metric cotype \(p \in (0, \infty)\) with constant \(\Gamma \in (0, \infty)\) if for every \(n \in \mathbb{N}\), there exists \(m = m(n, p, X) \in \mathbb{N}\) such that every function \(f: \mathbb{Z}_{2m}^n \to X\) satisfies
\begin{equation}\label{MCinq}
\left( \sum_{i=1}^n \sum_{x \in \mathbb{Z}_{2m}^n} d_X\left(f(x + m e_i), f(x)\right)^p \right)^{\frac{1}{p}} \leq \Gamma m \left( \frac{1}{2^n} \sum_{\varepsilon \in \{-1,1\}^n} \sum_{x \in \mathbb{Z}_{2m}^n} d_X\left(f(x + \varepsilon), f(x)\right)^p \right)^{\frac{1}{p}}.
\end{equation}
The space \((X, d_X)\) is said to have sharp metric cotype \(p \in (0, \infty)\) if there exist constants \(C, \Gamma \in (0, \infty)\) such that for every \(n \in \mathbb{N}\), there exists \(m \in \mathbb{N}\) with \(m \leq C n^{1/p}\) such that for every function \(f: \mathbb{Z}_{2m}^n \to X\), Equation~\eqref{MCinq} holds.
\end{definition}

The quantitative refinement from metric cotype to sharp metric cotype plays a crucial role in deriving non-embeddability results for spaces admitting sharp metric cotype. See \cite{MC} for further details.

Barycenter maps enable the formulation of martingales taking values in metric spaces. In \cite{CATCoarU}, such martingales were employed to show that \(\operatorname{CAT}(0)\) spaces possess sharp metric cotype \(2\). Nonlinear martingales are discussed, for instance, in \cite{NMinNPC}. The definition below is adapted from \cite{LipForBary}, where it is stated for probability measures with finite support, sufficient for establishing the sharp metric cotype \(2\) inequality in \(\operatorname{CAT}(0)\) spaces.

\begin{definition}[Martingale]
Let \((X, d_X)\) be a metric space that admits an \(\infty\)-barycenter map $\mathfrak{B}:\mathcal{P}_{\infty}\rightarrow X$. Let \(\Omega\) be a finite set and \(\mu: 2^{\Omega} \to [0,1]\) a probability measure with full support, meaning that \(\mu(\{\omega\}) > 0\) for every \(\omega \in \Omega\).
\begin{enumerate}
    \item Let \(\mathcal{F} \subseteq 2^\Omega\) be a \(\sigma\)-algebra. For each \(\omega \in \Omega\), let \(\mathcal{F}(\omega)\) denote the unique atom of \(\mathcal{F}\) containing \(\omega\). Given a function \(Z: \Omega \to X\), its \(\mu\)-conditional barycenter is the function \(\mathfrak{B}_\mu(Z \mid \mathcal{F}): \Omega \to X\) defined by
    \[
    \mathfrak{B}_\mu(Z \mid \mathcal{F})(\omega) = \mathfrak{B} \left( \frac{1}{\mu(\mathcal{F}(\omega))} \sum_{a \in \mathcal{F}(\omega)} \mu(a) \delta_{Z(a)} \right).
    \]
    
    \item Fix \(n \in \mathbb{N}\), and let \(\mathcal{F}_0 \subseteq \mathcal{F}_1 \subseteq \cdots \subseteq \mathcal{F}_n \subseteq 2^\Omega\) be a filtration. A sequence of functions \(\{Z_i: \Omega \to X\}_{i=0}^n\) is called a \(\mu\)-martingale with respect to the filtration \(\{\mathcal{F}_i\}_{i=0}^n\) if
    \[
    \mathfrak{B}_\mu(Z_i \mid \mathcal{F}_{i-1}) = Z_{i-1}, \quad \forall i \in \{1, \ldots, n\}.
    \]
\end{enumerate}
\end{definition}

The following inequality is a nonlinear analogue of a classical martingale inequality due to Pisier \cite{MartInUniB}, originally established for martingales taking values in uniformly convex Banach spaces. A generalization to metric spaces was given in \cite{LipForBary} for spaces \((X, d_X)\) that admit a barycenter map satisfying the Generalized Variance Inequality \eqref{VarInEq} for all \(z \in X\), with exponent \(p \in [1, \infty)\). Such spaces are referred to as \(p\)-barycentric in \cite{LipForBary}.

Moreover, a closer inspection of the proof in \cite{LipForBary} shows that if \eqref{VarInEq} holds only at a fixed point \(z \in X\), then the inequality stated below still holds at that same point. Therefore, by \cite[Lemma 2.1]{LipForBary} and Theorem~\ref{VarIn}, we obtain the following.

\begin{proposition}[Pisier's inequality]\label{PisIn}
    Let \(p \in [2, \infty)\), let \((X, d_X)\) be a uniquely geodesic space, and let \(z \in X\). Assume that the space \((X, d_X)\) is \(p\)-uniformly convex with respect to \(z\) with constant \(k > 0\), and that it admits an \(\infty\)-convex mean map \(\mathfrak{B} \colon \mathcal{P}_{\infty}(X) \to X\) which is invariant under midpoint contractions. Suppose \(\mu\) is a probability measure with full support on a finite set \(\Omega\), and let \(\{Z_i: \Omega \to X\}_{i=0}^n\) be a \(\mu\)-martingale with respect to a filtration \(\{\Omega, \emptyset\}=\mathcal{F}_0 \subseteq \cdots \subseteq \mathcal{F}_n \subseteq 2^\Omega\). Then, we have 
\begin{equation}\label{PisEq}
    \frac{2^{p-2}}{2^{p-1}-1} \frac{k}{2} \sum_{i=1}^n \int_{\Omega} d_X\left(Z_{i}, Z_{i-1}\right)^p d \mu \leqslant \int_{\Omega} d_X\left(Z_{n}, z\right)^p d \mu-\int_{\Omega} d_X\left(Z_{0}, z\right)^p d \mu
\end{equation}
In particular, if \((X, d_X)\) is a \(p\)-uniformly convex metric space with constant \(k > 0\) that admits a barycenter map as above, then Equation~\eqref{PisEq} holds for every \(z \in X\).
\end{proposition}
Combining Proposition~\ref{PisIn} with Property \eqref{conk}, we arrive at the following result.
\begin{corollary}
Let \((X, d_X)\) be the ball \(B_r(o)\) of radius \(0 < r < \pi / 2\), centered at the midpoint \(o\), in a complete \(\operatorname{CAT}(1)\) space. Suppose that \(\mu\) is a probability measure with full support on a finite set \(\Omega\), and let \(\{Z_i : \Omega \to X\}_{i=0}^n\) be a \(\mu\)-martingale with respect to a filtration \(\{\Omega, \emptyset\} = \mathcal{F}_0 \subseteq \cdots \subseteq \mathcal{F}_n \subseteq 2^\Omega\). Then, Equation~\eqref{PisEq} holds with \( z = o \), \( p = 2 \), and \( k = 2r \tan\left(\frac{\pi}{2} - r\right) \).
\end{corollary}
From Proposition~\ref{PisIn} and Property \eqref{conp}, we also derive the following.
\begin{corollary}
Let \((X, d_X)\) be a complete \(\operatorname{CAT}(0)\) space. Suppose that \(\mu\) is a probability measure with full support on a finite set \(\Omega\), and let \(\{Z_i : \Omega \to X\}_{i=0}^n\) be a \(\mu\)-martingale with respect to the filtration \(\{\Omega, \emptyset\} = \mathcal{F}_0 \subseteq \cdots \subseteq \mathcal{F}_n \subseteq 2^\Omega\). Then, Equation~\eqref{PisEq} holds for all \( z \in X \) with \( p = 2 \) and \( k = 2 \), or for \( p > 2 \) and \( k = \left( \frac{(p-1)^2}{10 p^2} \right)^{p-1} \).
\end{corollary}
Pisier's inequality was used in the proof of \cite[Theorem 5]{CATCoarU}, which implies that complete \(\operatorname{CAT}(0)\) spaces possess sharp metric cotype 2. From the same theorem, together with Theorem \ref{VarInCAT(0)}, we obtain the following generalization.

\begin{corollary}\label{MSforCAT}
Every complete \(\operatorname{CAT}(0)\) space has sharp metric cotype \(p\), for all \(p \in [2, \infty)\).
\end{corollary}

Since metric cotype is invariant under bi-Lipschitz embeddings, and every \(\operatorname{CAT}(1)\) space embeds bi-Lipschitzly into its Euclidean cone, which is a \(\operatorname{CAT}(0)\) space (see \cite[II.3.14]{MSofNPC}), the \(\operatorname{CAT}(1)\) version of Corollary~\ref{MSforCAT} follows formally from the \(\operatorname{CAT}(0)\) case.

\vspace{\baselineskip}

\noindent \textbf{Acknowledgments.}
This work is based on the author’s master’s thesis at the University of Bonn under the supervision of Professor Karl-Theodor Sturm, whose guidance is gratefully acknowledged.
\bibliography{references}

\end{document}